\documentclass[12pt]{article}
\newcounter{example} %example counter
\newtheorem{theorem}{Theorem}
\newtheorem{lemma}[theorem]{Lemma}

\renewcommand{\Re}{\mbox{Re\,}}

\newcommand{\ra}{\rightarrow}
\newcommand{\var}{\mbox{\rm var\,}}
\newcommand{\tr}{\mbox{\rm tr\,}}
\newcommand{\F}{{\cal F}}
\newcommand{\M}{{\cal M}}
\renewcommand{\H}{{\bf H}}
\newcommand{\R}{{\bf R}}

\newcommand{\be}{\begin{equation}}
\newcommand{\ee}{\end{equation}}
\newcommand{\xx}{x_1,\ldots,x_N}
\newcommand{\dxx}{dx_1 \cdots dx_N}

\begin{document}
	\date{}

\title{Distribution Functions for Random Variables for Ensembles of 
Positive Hermitian Matrices}

\author{Estelle L. Basor \thanks{ebasor@calpoly.edu. 
          Supported in part by NSF Grant DMS-9623278. }   \\
               Department of Mathematics\\
               California Polytechnic State University\\
               San Luis Obispo, CA 93407, USA}

\maketitle

\begin{abstract}
Distribution functions for random variables that depend on a parameter
are computed asymptotically 
for ensembles of positive Hermitian matrices. The inverse Fourier transform
of the distribution is shown to be a Fredholm determinant of a certain 
operator that is an analogue of a Wiener-Hopf operator. The asymptotic 
formula shows that up to the terms of order $o(1)$, the distributions are Gaussian.
\end{abstract}

\section{Introduction}
In the theory of random matrices one is led naturally to consider the 
probability distribution on the set of eigenvalues of the matrices. For 
$N \times N$ random Hermitian matrices one can show that under 
reasonable assumptions, the probability density that the eigenvalues 
$\lambda_1, \ldots, \lambda_N$ lie in the intervals
\[ (x_1, x_1 + dx_1), \ldots, (x_N, x_N + dx_N) \]
is given by the formula

\begin{equation}
  P_N(x_1,\ldots,x_N) = \frac{1}{N!}\det K(x_i,x_j) \left. \right|_{i,j=1}^N
\end{equation} \label{eqn1.1}

where

\begin{equation}
  K_N(x,y)= \sum_{i=0}^{N-1}\phi_i(x) \phi_i(y)
\end{equation} 

and  $\phi_i$ is obtained by orthonormalizing the sequence 
$\left\{x^ie^{-x^2/2} \right\}$ over $\R$.

For $N \times N$ positive Hermitian matrices the probability density has 
the same form except that $\phi_i$ is replaced by the functions obtained 
by orthonormalizing the sequence $\left\{ x^{\nu / 2} e^{-x/2} x^i \right\} $ over
$\R ^+ $. We will not describe here exactly how these particular 
densities arise but instead refer the reader to \cite{M}.

We can define a random variable on the space of eigenvalues by considering
$f(x_1, \ldots, x_N)$ where $f$ is any symmetric function of the $x_i$'s. 
A particular case of interest is a random variable of the form
$\sum_{i=1}^{N}f(x_i)$, where $f$ is a function of a real variable. Such 
a random variable is generally called a linear statistic.

In previous work \cite{M,B,BT,Johansson}, the variance of the random 
variable was computed in the large $N$ limit. More precisely, the 
function $f$ and the kernel $K_N(x,y)$ were suitably rescaled so that 
the limit as $N \rightarrow \infty$ of the variance could be computed.
The precise details of this are in the next section.

Our goal in this paper is to compute the distribution function for a 
class of the linear statistics that depend on a parameter $\alpha$.
We now describe the sections of the paper and main results. In the next 
section we outline the random matrix theory and show how the distribution 
functions can be computed using Fredholm determinants. In section 3 we 
replace the function $f(x)$ in the linear statistic by 
$f_\alpha (x) = f(x/\alpha)$. For random variables of this type we show 
that the inverse Fourier transform of the distribution function 
$\check{\phi}(k)$
 has an 
asymptotic expansion of the form
\be
  \check{\phi}(k) \sim e^{ak^2 + bk}
\ee
as $ \alpha \rightarrow \infty $. This of course implies that the actual 
distribution is asymptotically Gaussian.
Here $a$ and $b$ depend on $f$ and $\alpha$. This is proved for both the 
Hermitian matrices and positive Hermitian matrices. In the latter case 
with $\nu =-1/2$, a very simple proof is given in Section 3. For $\nu > 
-1/2$, a completely different proof is obtained in Section 4.
 
  Most of the results are obtained by 
using simple operator theory identities in the theory of Wiener-Hopf 
operators. The central idea is that the various quantities which yield 
information about random variables can all be computed in terms of traces 
or determinants of integral operators. Some of the computations lead 
directly to a familiar problem in the theory of Wiener-Hopf operators, 
while others require modifications and generalizations of these results.
%end page 1
\section{Preliminaries}
In this section we show how to compute the mean, variance, and inverse Fourier 
transform of the distribution of the random variable. Computations for 
the mean and variance have been given before in many places. However, we 
reproduce all of these here for completeness sake and also to highlight 
the use of operator theory ideas.

We begin by considering $P_N$ for $N \times N$ random Hermitian matrices. 
We want to consider large matrices and thus we let $N \rightarrow \infty$,
but this leads to a trivial result unless we rescale $K_N$ in a 
particular way. We replace $K_N(x,y)$ with
\be \frac{1}{\sqrt{2N}} 
K_N\left(\frac{x}{\sqrt{2N}},\frac{y}{\sqrt{2N}}\right). \label{eqn2.1}
\ee  
Rescaling $K_N$ is equivalent to rescaling the mean spacing of the 
eigenvalues. (See \cite{TW} for details.) 

From the theory of Hermite polynomials it is 
easy to see that as $N \ra \infty$,
\be
 \frac{1}{\sqrt{2N}}K_N\left(\frac{x}{\sqrt{2N}},\frac{y}{\sqrt{2N}}\right)
 \ra \frac{\sin(x-y)}{\pi (x-y)}.
\ee
This last function is known as the sine kernel. Now consider a random 
variable of the form
\[ \sum_{i=1}^{N}f(x_i\sqrt{2N}) \]
where in all that follows $f$ is a continuous real-valued function 
belonging to $L_1(\R)$ and which vanishes at $\pm \infty$.
The appearance of the $\sqrt{2N}$ should not be surprising here since the 
above rescaling spreads out the eigenvalues and hence should be reflected 
in the random variable. The mean $\mu_N$ is
\be
  \int \cdots \int \sum_{i=1}^{N}f(x_i\sqrt{2N}) P_N(x_1,\ldots,x_N)
  dx_1 \cdots dx_N .  \label{eqn2.2}
\ee
Now the function $P_N$ has the important property \cite{M}
\be
  \frac{N!}{(N-n)!}\int \cdots \int P_N(x_1,\ldots,x_n,x_{n+1},
  \ldots, x_N) dx_{n+1} \cdots dx_N = \det K(x_i,x_j) \left.
  \right|_{i,j=1}^n . \label{eqn2.3}
\ee
Thus, (\ref{eqn2.2}) is easily seen to be
\be
  \int_{-\infty}^{\infty}f(x\sqrt{2N})K_N(x,x)\, dx
  \label{eqn2.4}
\ee  %page 2
which, after changing $x$ to $x/\sqrt{2N}$, becomes
\[
  \int_{-\infty}^{\infty}f(x)\frac{1}{\sqrt{2N}} 
  K_N\left(\frac{x}{\sqrt{2N}},\frac{y}{\sqrt{2N}}\right) \, dx.
\]
Thus, as $N \ra \infty$, 
\be
  \mu_N \ra \mu = \int_{-\infty}^{\infty}f(x)K(x,x) \, dx
\ee
where $K(x,y)$ is the sine kernel.

A very similar computation for the variance $\var\!_N f$, again using (\ref{eqn2.3}),
yields
\be
  \var f:= \lim_{N \ra \infty} \var\!_N f =
  -\int \int f(x)f(y)K^2(x,y) \, dx \, dy +
  \int f^2(x) K(x,x) \, dx .
  \label{eqn2.6}
\ee
Both the mean and the variance can be interpreted as traces of certain 
Wiener-Hopf operators. To see this, consider the operator $A(f)$ on
$L_2(-1,1)$ with kernel
\be
  \frac{1}{2\pi}\int_{-\infty}^{\infty}f(t) e^{-it(x-y)} \, dt.
  \label{eqn2.7}
\ee
This operator can easily be seen to be the product 
$\F  \M_f \F^{-1} P$ where $P  g = \chi_{(-1,1)} g$, 
$\M_f g  = fg$ and $\F$ is the Fourier transform.
A moment's thought shows that $\mu = \tr\{A(f) \}$ and
$\var f = \tr \{ A(f^2) - (A(f))^2 \}$.

A more difficult, yet also straightforward problem, is to find an 
expression for the distribution function of a random variable of this 
type. A fundamental formula from probability theory shows that
if we call the probability distribution function $\phi_N$, then
\be
  \check{\phi}_N(k) = 
  \int_{-\infty}^{\infty}\cdots \int_{-\infty}^{\infty}
  e^{ik\sum_{j=1}^{N}f(x_j \sqrt{2N})}
  P_N(x_1,\ldots,x_N) dx_1 \cdots dx_N .
  \label{eqn2.8}
\ee
Thus,
\begin{eqnarray*}
	\check{\phi}_N(k) & = & \int_{-\infty}^{\infty}\cdots\int_{-\infty}^{\infty}
	\prod_{j=1}^{N}e^{ikf(x_j\sqrt{2N})}P_N(\xx)\,\dxx  \\
   & = & \int_{-\infty}^{\infty}\cdots\int_{-\infty}^{\infty}
 	 \prod_{j=1}^{N}((e^{ikf(x_j\sqrt{2N})}-1)+1)P_N(\xx)\,\dxx  \\
   & = & \int_{-\infty}^{\infty}\cdots\int_{-\infty}^{\infty}
	 \{1 + \sum_{j=1}^{N}(e^{ikf(x_j\sqrt{2N})}-1)
	 +\sum_{ j < l }^{N}
	 (e^{ikf(x_j\sqrt{2N})}-1)(e^{ikf(x_l\sqrt{2N})}-1)+\ldots \} \\
   &  & \times P_N (\xx)\,\dxx  \\
   & = & 1+\frac{1}{1!}\int_{-\infty}^{\infty}(e^{ikf(x\sqrt{2N})}-1)
	 K_N(x,x)\, dx   \\
   &  &	\mbox{} + \frac{1}{2!}\int_{-\infty}^{\infty}\int_{-\infty}^{\infty}
 	 (e^{ikf(x_1\sqrt{2N})}-1)(e^{ikf(x_2\sqrt{2N})}-1)
	 \det ( K_N(x_j,x_l)) \left. \right|_{1 \leq j,l \leq 2} \, dx_1 \, dx_2
	  \\
   &  &\mbox{} + \cdots + \frac{1}{N!}\int_{-\infty}^{\infty}\cdots\int_{-\infty}^{\infty}
	 \prod_{j=1}^{N}(e^{ikf(x_j\sqrt{2N})}-1) P_N(\xx)\, \dxx .	 
\end{eqnarray*}	%end page3
In each integral we rescale to obtain
\begin{eqnarray}
	\check{\phi}_N(k) & = & 1 + \frac{1}{1!}\int_{-\infty}^{\infty}K'(x_1,x_1) \, dx_1 
  + \frac{1}{2!}\int_{-\infty}^{\infty}\int_{-\infty}^{\infty}
  K'(x_1,x_2)\, dx_1 \, dx_2    \nonumber \\
   & & + \cdots +
  \frac{1}{N!}\int_{-\infty}^{\infty}\cdots\int_{-\infty}^{\infty}
  K'(\xx) \, \dxx
\end{eqnarray}
where 
\be
  K'(x_1,\ldots,x_n) = \det \left( 
  (e^{ikf(x_j)}-1)K_N(\frac{x_j}{\sqrt{2N}},\frac{x_l}{\sqrt{2N}})
  \frac{1}{\sqrt{2N}}\right)_{1 \leq j,l \leq n} .
\ee
Letting $N \ra \infty$ we see this is the formula for the Fredholm 
determinant  $\det (I + K)$ where $K$ has kernel
\be
  K(x,y) = (e^{ikf(x)}-1)\frac{\sin(x-y)}{\pi(x-y)}.
\ee
As before we can express this last quantity in terms of the operator 
$A(\sigma)$
\be
  \check{\phi}(k)=\lim_{N \ra \infty} \check{\phi}_N (k) = 
  \det (I + A(\sigma))
\ee
where $\sigma(x) = e^{ikf(x)}-1 $.
	
The preceding computations can all be carried out in the case of positive 
Hermitian matrices. In this case we replace $K_N(x,y)$ with
\[ \frac{1}{4N}K_N(\frac{x}{4N},\frac{y}{4N})  \]
and from the theory of Laguerre polynomials we see that as $N \ra \infty$
\be
  \frac{1}{4N}K_N(\frac{x}{4N},\frac{y}{4N}) \ra 
  \frac{{J_\nu}(\sqrt{x})\sqrt{y}{J_\nu}'(\sqrt{y})-\sqrt{x}{J_\nu}'(\sqrt{x})
    J_\nu(\sqrt{y})}{2(x-y)} \label{eqn2.10}
\ee
where $J_\nu$ is the Bessel function of order $\nu$. The details of this 
are found in \cite{TW3}. The rescaling here forces the eigenvalue density 
to be bounded near zero and is called ``scaling at the hard edge.'' The 
kernel (\ref{eqn2.10}) is known as the Bessel kernel.

We can again write the mean, the variance, and the Fourier transform of 
the distribution in terms of operators. This time the relevant operator 
$B(f)$ is defined on $L_2(0,1)$ with kernel given by
\be
  K(x,y)= \int_{0}^{\infty}t\sqrt{xy}f(t)J_\nu(tx)J_\nu (ty) \, dt.
  \label{eqn2.11}
\ee
If we begin with the linear statistic (the $\sqrt{x}$ is merely for 
convenience, and we again assume that $f$ is continuous, in $L_{1}(\R ^{+} )$
and vanishes at $+ \infty$)
\be
  \sum_{i=1}^{N}f(\sqrt{x_i 4N}),
\ee
then nearly identical computations show that
\begin{eqnarray*}
 \mu & = & \tr B(f) \\
 \var f & = & \tr\{ B(f^2) - (B(f))^2 \} \\
 \check{\phi}(k) & = & \det (I + B(\sigma))
\end{eqnarray*}
where $\sigma = e^{ikf(x)} -1 $. 
We summarize these results in the following:  %end page4
\begin{theorem}
(a) Given a random variable of the form $\sum_{i=1}^{N}f(x_i\sqrt{2N})$
defined on the space of eigenvalues of $N \times N$ Hermitian matrices 
with probability distribution given in (\ref{eqn1.1}), we have
\[\begin{array}{rcccl}
	\mu & := & \lim_{N \ra \infty} \mu_N & = & \tr(A(f))  \\
	\var f & := & \lim_{N \ra \infty} \var \!_N f & = & \tr \{ 
	A(f^2)-(A(f))^2\}  \\
	\check{\phi}(k) & := & \lim_{N \ra \infty} \check{\phi}_N (k) & = & 
	\det(I + A(\sigma))
\end{array} \]
where 
$ \sigma(x) = e^{ikf(x)} -1  $.

(b) Given a random variable of the form $\sum_{i=1}^{N}f(\sqrt{x_i4N})$
defined on the space of eigenvalues of positive $N \times N$ Hermitian 
matrices, we have
\[\begin{array}{rcccl}
	\mu & := & \lim_{N \ra \infty} \mu_N & = & \tr(B(f))  \\
	\var f & := & \lim_{N \ra \infty} \var \!_N f & = & \tr \{ 
	B(f^2)-(B(f))^2\}  \\
	\check{\phi}(k) & := & \lim_{N \ra \infty} \check{\phi}_N (k) & = & 
	\det(I + B(\sigma))
\end{array} \]
where 
$ \sigma(x) = e^{ikf(x)} -1 $.
\end{theorem}
 
When linear statistics are considered \cite{B,SMMP}, one is often 
concerned with a statistic of the form $\sum_{i=1}^{N}f(x_i/\alpha)$
where $\alpha$ is a real parameter approaching infinity. This is the case,
for example, 
in the study of disordered conductors where large $\alpha$ corresponds 
to a high density metallic regime. The above formulas still hold, of 
course, but now they depend on the parameter. We will call the operators 
that depend on the parameter $\alpha$ by $A_\alpha(f)$ and $B_\alpha(f)$,
respectively. In the next sections we will compute the mean, variance, 
and distribution function asymptotically as $\alpha \ra \infty$.

\section{The Mean, Variance, and Distribution Function as $\alpha \ra 
\infty$}

For random Hermitian matrices, computing the various limits are 
applications of the continuous analogues of the Strong Szeg\"o Limit 
Theorem. For then, $A_\alpha(f)$ is just the classical Wiener-Hopf 
operator defined on the interval $(-\alpha,\alpha)$, and all of the 
quantities are known asymptotically as $\alpha \ra \infty$. We provide 
the answers here for completeness.

\begin{theorem}
 Assume that $f \in L_1(\R)$ is continuous, and vanishes at $\pm \infty$
  and that in addition its  Fourier transform  $\hat{f}$ satisfies
 
 \[  \int_{-\infty}^{\infty}|x| |\hat{f}(x)|^2 \, dx \; < \; \infty.  
 \]
 
 Then
 \begin{eqnarray*}
  \mu & = & \frac{\alpha}{2\pi}\int_{-\infty}^{\infty}f(x) \, dx \\
  \var f & = & 2 \int_{0}^{\infty}x \hat{f}(x)\hat{f}(-x)\, dx 
  \,+\,o(1)
 \end{eqnarray*}
 and 
 \[ 
    \check{\phi}(k) \sim \exp \left\{ \frac{\alpha}{2\pi}\int_{-\infty}^{\infty}
    ikf(x) \, dx  - k^2 \int_{0}^{\infty}x \hat{f}(x)\hat{f}(-x)\, dx 
    \right\}.
 \] 
  \end{theorem}
%end page 5  
The Bessel case is significantly more complicated. There is no 
corresponding Szeg\"o type theorem. We begin by computing the mean. The 
operator $B_\alpha(\sigma)$ has kernel
\[
  \int_{0}^{\infty}\sqrt{xy}tf(t/\alpha)J_\nu(tx)J_\nu(ty) \, dt .
\]
Thus the mean $\mu$ is given by
\begin{eqnarray}
	\mu & = & \int_{0}^{1}\int_{0}^{\infty}xtf(t/\alpha)J_\nu^2(tx)
	           \, dt \, dx  \nonumber \\
	 & = & \alpha^2\int_{0}^{\infty}\int_{0}^{1}xtf(t)J_\nu^2(\alpha tx)
	         \, dx \, dt  \nonumber \\
	 & = & \alpha^2\int_{0}^{\infty}f(t)\int_{0}^{1}xt J_\nu^2(\alpha t x)
	       \, dx \, dt . \label{eqn3.1}
\end{eqnarray}
Now
\[ 
  \int_{0}^{1}x J_\nu^2(\alpha t x) \, dx = 
  \frac{1}{2} \left\{ J_\nu ^2(\alpha t) - 
  J_{\nu + 1}(\alpha t) J_{\nu - 1}(\alpha t)  \right\}
\]
and
\[
  J_{\nu -1} (\alpha t) = -J_{\nu + 1} (\alpha t) 
    + \frac{2\nu}{\alpha t}J_\nu(\alpha t) .
\]

Therefore the integral (\ref{eqn3.1}) becomes
\[
  \alpha \int_{0}^{\infty} f(t) \frac{\alpha t}{2} \left\{ J_\nu^2(\alpha t)
  + J_{\nu + 1}^2(\alpha t) - \frac{2\nu}{\alpha t} 
  J_{\nu+1}(\alpha t) J_\nu(\alpha t) \right\} \, dt
\]
or
\be
  \alpha \int_{0}^{\infty} f(t) \frac{\alpha t}{2} \left\{ J_\nu^2(\alpha t)
  + J_{\nu + 1}^2(\alpha t) \right\} \, dt - \alpha \nu \int_{0}^{\infty}
  f(t) J_{\nu + 1}(\alpha t)J_\nu(\alpha t) \, dt.
\ee

The first integral equals 
\be
 \frac{\alpha}{\pi} \int_{0}^{\infty} f(t)\,dt \,\,+o(1)
 \ee
 which can be easily seen by using the asymptotic properties of Bessel 
 functions. The second integral is asymptotically 
 \[
  \frac{\nu}{2}f(0) +o(1).
\]
This uses the identity $\int_{0}^{\infty}J_{\nu +1}(x)J_{\nu }(x)\,dx 
=\frac{ 1}{2}.$

Thus we have

\be
 \mu = \frac{\alpha}{\pi}\int_{0}^{\infty}f(t)\,dt -\frac{\nu}{2}f(0) +o(1).
 \ee
 
 For the variance we refer to \cite{BT} where the calculation was already 
 done. There it was shown that
 \be
 \mbox{var}f \sim \frac{1}{\pi^{2}} \int_{-\infty}^{\infty} |M(f)(2iy)|^{2}y 
  \tanh (\pi y) dy.
  \ee
  %end page 6
  We note however, that this can also be written as 
  \be
   \mbox{var}f \sim \frac{1}{\pi^{2}} \int_{0}^{\infty}x (C(f)^{2})\,dx
  \ee
  where $C(f)(x) = \int_{0}^{\infty}f(y)\cos (xy)dy $ denotes the cosine transform of $f$. This is an exercise 
  involving the properties of the Mellin transform, and we leave it to the 
  reader. 
  
  To compute the distribution function, we first turn our attention to 
  the case where $\nu = -1/2$. Our operator $B_{\alpha}(\sigma)$ has 
  kernel

  \begin{eqnarray*} 
    &  &  \frac{2}{\pi} \int_{0}^{\infty} \sigma(t/\alpha) \cos xt \cos 
    yt dt    \\
  & = & \frac {1}{\pi}\int_{0}^{\infty} \sigma(t/\alpha) 
  (\cos((x-y)t)+\cos((x+y)t))\,dt\\
 & = & \frac{\alpha}{\pi}(C(\sigma)((x-y)\alpha)+ C(\sigma)((x+y)\alpha)).
  \end{eqnarray*}

  This is unitarily equivalent to the operator on $L_{2}(0,\alpha)$ with 
  kernel
  \be
    \frac{1}{\pi}(C(\sigma)(x-y) + C(\sigma)(x+y) ).
    \ee
    
    The operator with kernel $\frac{1}{\pi}(C(\sigma)(x-y)) $ is the 
    finite Wiener-Hopf operator, usually denoted as $W_{\alpha}(\sigma)$, and the 
    operator with kernel $\frac{1}{\pi}(C(\sigma)(x+y)) $ is the Hankel 
    operator $H_{\alpha}(\sigma)$. (The only difference between this 
    definition of a finite Wiener-Hopf operator and the one given earlier for 
    $A_{\alpha}$ is the difference in the domain. The two are unitarily 
    equivalent.) If we consider the operators 
    on $L_{2}(0,\infty)$ in what follows, we will denote them by 
    $W(\sigma)$ and $H(\sigma)$ respectively. Also, whenever it is 
    necessary to consider the extension of
    $\sigma$ to the entire real axis, it will always be the even extension.
    
    Thus the problem of finding the distribution function asymptotically 
    becomes the same as computing the Fredholm determinant 
    $\det(I+B_{\alpha}(\sigma))=\det(I+W_{\alpha}(\sigma)+H_{\alpha}(\sigma))$ 
    asymptotically. 
    To do this we need some basic facts about Wiener-Hopf operators and 
    we collect them in the following theorem. These are well-known and 
    can all be found in \cite{BS}.
    
    \begin{theorem}
    a) Suppose $\phi$ and $\psi$ are even bounded functions in $L_{1}(\R)$. Then
 $$W(\phi)H(\psi) + H(\phi)W(\psi) = H(\phi \psi)$$
 \center{and}
  $$W(\phi)W(\psi) = W(\phi \psi) - H(\phi)H(\psi).$$
 b) Suppose $\phi$ and $\psi$ are bounded functions in $L_{1}(\R)$.
     If the Fourier transform $\hat{\phi}(x)$ vanishes for $x$ negative, 
    then $ W(\psi)W(\phi) =W(\phi \psi)$ and
     if $\hat{\phi}(x)$ vanishes for $x$ positive, then $ 
    W(\phi)W(\psi) =W(\phi \psi)$.
    \end{theorem}
    
    We define $W(\sigma)$ and $H(\sigma)$ with $\sigma = 1+f$ and $f$ in $L_{1}$ by 
    $W(\sigma) = I+W(f)$ and $H(\sigma) = H(f)$. Both of these 
    definitions are natural when thought of in a distributional setting, 
    and the above theorem holds with these definitions as well.
    
    The next theorem is of primary importance in the computations that 
    follow.
    
    \begin{theorem}
    Suppose $\phi = 1+f, \phi^{-1} = 1+g$ where $f $ and $g$ are bounded even 
    functions. Then the inverse of $W(\phi) +H(\phi)$ is $W(\phi^{-1}) + 
    H(\phi^{-1})$.
    \end{theorem}
    
    Proof: Using Theorem 3 parts a) and b) we have, 
    \begin{eqnarray*} 
     &   & (W(\phi) +H(\phi)) (W(\phi^{-1}) + H(\phi^{-1}))\\
     & = & W(\phi)W(\phi^{-1}) + H(\phi)W(\phi^{-1}) 
     +W(\phi)H(\phi^{-1})+H(\phi)H(\phi^{-1})\\
     & = &I-H(\phi)H(\phi^{-1}) + H(\phi \phi^{-1}) +H(\phi)H(\phi^{-1})\\
     & = & I + H(1) = I.\\
     \end{eqnarray*} %end page 7
    The same computation holds for $(W(\phi^{-1}) + H(\phi^{-1}))(W(\phi) 
     +H(\phi))$, and so we have shown that these operators are inverses of 
     each other.
    
     It is well known from the theory of Wiener-Hopf operators that under 
     appropriate conditions $\det(I+W_{\alpha}(\sigma))$ has the 
     asymptotic expansion $G(\sigma)^{\alpha}E(\sigma)$ where $$G(\sigma) 
     =\exp \frac{1}{2 \pi}\int_{-\infty}^{\infty} \log (1+\sigma(\xi) ) d 
     \xi$$ 
     and $E(\sigma) = \det (W(\phi)W(\phi^{-1}))$ with $\phi = 
     1+\sigma.$ This is simply another version of Theorem 2. With 
     additional assumptions on $\phi$, it is very easy to adapt this proof 
     to the Bessel case $\nu = -1/2$ to show that 
     \be
      \det(I+W_{\alpha}(\sigma) + H_{\alpha}(\sigma)) \sim 
     G(\sigma)^{\alpha}E'(\sigma) \label{A}
     \ee
     and $E'(\sigma) = 
     \det((W(\phi)+H(\phi))W(\phi^{-1}))$. Thus to compute the 
     distribution, we need to know the form of the above determinant. This 
     is contained in the next theorem.
     
     \begin{theorem}
     Suppose $\sigma = e^{ikf}-1$ where $f$ is even, continuous, piecewise 
     $C^{2}$ and vanishes at infinity. Suppose also that $f\in L_{1}$ and the 
     function $$ \xi \ra (1+\xi^{2})(|f''(\xi)|+|f'(\xi)|^{2}) \in 
     L_{2}.$$ Then as $\alpha \ra \infty$, we have 
      
      \be
         \det(I+W_{\alpha}(\sigma) + H_{\alpha}(\sigma)) \sim \exp \{ 
         \frac{\alpha}{ \pi} \int_{0}^{\infty} ikf(x)\,dx 
         +\frac {ik}{4} f(0) -\frac{k^{2}}{2 \pi^{2}}\int_{0}^{\infty} 
         x |C(f)(x)|^{2}\,dx \}.  
      \ee
      \end{theorem}
      Proof: The conditions on $\sigma$ ensure that the above integrals converge, 
      and that the operators $H(\phi)$ and $H(f)$ are trace class. The 
      reader is referred to \cite{BW} for details. These assumptions also 
      guarantee that (27) holds. It is also easy to see that  
      $G(\phi) = \exp \{ \frac{\alpha} { \pi} 
      \int_{0}^{\infty}ikf(x)\,dx\}.$ 
      To complete the proof we need a concrete representation 
      for $ \det((W(\phi)+H(\phi))W(\phi^{-1})).$ Define $$h(k) = \log 
      \det ((W(\phi)+H(\phi))W(\phi^{-1}))$$ where $\phi = e^{ikf}.$  Let 
      $h(k) = \log \det ((W(\phi)+H(\phi))W(\phi_{-1})).$ We need to show 
      the second derivative of $h$ is constant in $k$. A standard formula
      \cite{GK} yields 
      
      \begin{eqnarray*}
         h'(k) & = & \mbox{tr}((W(\phi^{-1}))^{-1}(W(\phi)+H(\phi))^{-1} \times
      \frac{d(W(\phi)+H(\phi))W(\phi^{-1})}{dk})\\
      &  = & \mbox{tr} ((W(\phi^{-1}))^{-1}(W(\phi)+H(\phi))^{-1})\\
      &    &  \times \{(W(\phi)+H(\phi))W(\phi^{-1}(-if))+W(\phi if)W(\phi^{-1})+H(\phi if)W(\phi^{-1}) \}\\
      & = & 
      \mbox{tr}\{(W(\phi^{-1}))^{-1}W(\phi^{-1}(-if))
      +(W(\phi^{-1}))^{-1}W(\phi^{-1})W(\phi if)W(\phi^{-1})\\
      &   &+ (W(\phi^{-1}))^{-1}W(\phi^{-1})H(\phi if)W(\phi^{-1}) + 
      (W(\phi^{-1}))^{-1}H(\phi^{-1}) W(\phi if)W(\phi^{-1})\\
      &   & +(W(\phi^{-1}))^{-1}H(\phi ^{-1})H(\phi if)W(\phi^{-1}) \}.\\
     \end{eqnarray*} 
     This uses Theorem 4. Simplifying further  and using the fact 
      that $H(\phi^{-1})$ is trace class we have
      \begin{eqnarray*} 
      h'(k) & = &\mbox{tr}\{ (W(\phi^{-1}))^{-1}W(\phi^{-1}(-if))+W(\phi 
      if)W(\phi^{-1})\\
       &   & + H(\phi if)W(\phi^{-1})+H(\phi^{-1})W(\phi 
      if)+H(\phi^{-1})H(\phi if)\}.
      \end{eqnarray*} %end page 8
    Now apply Theorem 3, part a) and the fact that 
      $\mbox{tr}(AB)=\mbox{tr}(BA)$ to find
      \begin{eqnarray*}
      h''(k) & = & \mbox{tr} \{(W(\phi^{-1}))^{-1}W((\phi^{-1})(if)^{2})\\
      &  & -(W(\phi^{-1}))^{-1}W(\phi^{-1} 
      (-if))(W(\phi^{-1}))^{-1}W(\phi^{-1}(-if))\}.
      \end{eqnarray*}
      The conditions on $\phi$ guarantee that the function $\phi $   
        has a factorization $\phi  = (g_{-}+1)(g_{+}+1)$ such that the 
        Fourier transforms of $g_{+}$ and $g_{-}$ vanish for positive 
        and negative real values respectively. Then using Theorem 3, 
        part b),it is easy to see that we can write
        $$ W(\phi) = W(g_{-}+1)W(g_{+}+1), 
        W(\phi^{-1})^{-1}=W(g_{+}+1)W(g_{-}+1).$$
        A repeated application of these identities allows us to write 
        $h''(k) = \mbox{tr}H(if)H(if),$ and $h''(k)$ is independent 
        of $k$. Thus at this point we have $h(k) 
        = ak^{2} +bk + c$ where $2a = -\mbox{tr}((H(f))^2.$ A direct 
        computation shows that $a = - \frac{1}{2\pi^{2}} \int_{0}^{\infty} 
        x|C(f)(x)|^2\,dx.$ To compute $b$, notice that $h'(0)$ is 
        $\mbox{tr} H(if) = \frac{i}{2 \pi}\int_{0}^{\infty}C(f(x))\,dx.$
        Also $h(0) = \mbox{tr} \log(I) = 0.$ Thus the last theorem holds.
        
   \section{The General Case}
   
    In this section we show that under certain conditions, the 
   distribution function for general $\nu$ has the same form as in the case
    of $\nu = -1/2.$ The only difference is in the mean which was 
    computed in the last section. The attack on the problem is entirely 
    different here. Instead of computing determinants asymptotically,
     we compute the traces of the 
    operators
    $(B_{\alpha}(\sigma))^{n}$ and then piece together the answers to get 
    an answer for the trace of $\log(I + B_{\alpha}(\sigma))$ and from that to the
    desired determinant.

   To begin  we need to show that $\tr 
   f(B_{\alpha}(\sigma))$ makes sense for a class of analytic functions 
   $f$. Just as we can associate the Wiener-Hopf operator with the Fourier 
   transform and a multiplication operator, we can also write 
   $$B_{\alpha}(\sigma) = P \H \M_{\sigma} \H  $$
   where $\H$ is the Hankel transform and $P$ is the projection on 
   $L_{2}(0,1)$. Since the Hankel transform is unitary on $L_{2}(0,\infty)$
   (\cite{UN}), the operator norm $||B_{\alpha}(\sigma)||$ is less than the 
   infinity norm $||\sigma ||_{\infty}$ of $\sigma$. Thus 
   $f(B_{\alpha}(\sigma))$ is defined for $f$ analytic on a disk centered 
   at the origin with radius $||\sigma||_{\infty}+\delta , \delta >0.$ 
   The operator $B_{\alpha}(\sigma)$ is also trace class for $\sigma$ in 
   $L_{1}$ by Mercer's Theorem (\cite{GK} Ch.III) as is 
   $f(B_{\alpha}(\sigma))$ for $f$ satisfying the above and $f(1)=0$.

We need  some lemmas that will prove to be useful. These may be 
    known already, but we include them for completeness.
    
    \begin{lemma}
    Suppose $-1 <  p < 1, \,\, 0< \lambda, \delta <1, \mu< 0, p + \mu + 
    \delta < 0$ and $0 < t <1.$ Then 
    \be
    \int_{0}^{\infty} s^p (1+s)^{\mu} |1-s|^{-1+ 
    \lambda}|1-ts|^{-1+\delta}\,\,ds \leq A \max (|1-t|^{-1+\lambda}, 
    |1-t|^{-1+\delta}) \max (t^{-\lambda},t^{-p-\lambda})
   \ee 
   where $A$ is some constant independent of $t$.
   \end{lemma}
   
   Proof: We have 
   \begin{eqnarray*}
   &  & \int_{0}^{\infty} s^p (1+s)^{\mu} |1-s|^{-1+ 
    \lambda}|1-ts|^{-1+\delta}\,\,ds\\
    & = & t^{-1+\delta}\int_{0}^{1} s^p (1+s)^{\mu} |1-s|^{-1+ 
    \lambda}|1/t - s|^{-1+\delta}\,\,ds\\
    & & +t^{-1+\delta}\int_{1}^{1/t} s^p (1+s)^{\mu} |1-s|^{-1+ 
    \lambda}|1/t-s|^{-1+\delta}\,\,ds\\
    & & + t^{-1+\delta}\int_{1/t}^{\infty} s^p (1+s)^{\mu} |1-s|^{-1+ 
    \lambda}|1/t-s|^{-1+\delta}\,\,ds.
    \end{eqnarray*}
   We consider each of the above integrals. In each, $A$ is a possibly 
    different constant independent of $t$ but can depend on the other 
    parameters.
    %end page 9
    First, 
    \begin{eqnarray*}
    &  &\int_{0}^{1} s^p (1+s)^{\mu} |1-s|^{-1+ 
    \lambda}|1/t-s|^{-1+\delta}\,\,ds\\
    & \leq & |1/t -1|^{-1+\delta}\int_{0}^{1}s^p(1+s)^{\mu}|1-s|^{-1+\lambda}ds\\
    & \leq & A|t-1|^{-1+\delta}t^{1-\delta}.
   \end{eqnarray*}
   Next,
   \begin{eqnarray*}
   & & \int_{1}^{1/t} s^p (1+s)^{\mu} |1-s|^{-1+ 
    \lambda}|1/t-s|^{-1+\delta}\,\,ds \\
    & \leq & A\max (1,t^{-p})\int_{1}^{1/t}|1-s|^{-1+\lambda}|1/t 
    -s|^{-1+\delta} ds\\
    & = & A\max (1, t^{-p})|1/t-1|^{-1+\lambda + \delta}\\
    & = & A\max (1, t^{-p})t^{-\lambda - \delta +1}|1-t|^{-1+\lambda + 
    \delta}\\
    & \leq &  A\max (1, t^{-p})t^{-\lambda - \delta +1}|1-t|^{-1+\lambda}.\\
   \end{eqnarray*}
    Finally, 
   \begin{eqnarray*}
   &  & \int_{1/t}^{\infty} s^p (1+s)^{\mu} |1-s|^{-1+ 
    \lambda}|1/t-s|^{-1+\delta}\,\,ds.\\
    & \leq & |1-1/t|^{-1+\lambda}\int_{1/t}^{\infty} s^{p+\mu}|1/t-s|^{-1+\delta}
    \,\,ds\\
    &  \leq & |1-1/t|^{-1 + \lambda}t^{-p-\mu -\delta}A\\
    & =  & A |t-1|^{-1+\lambda}t^{-p-\mu -\delta -\lambda +1}.
    \end{eqnarray*}
 Putting this together we have that the original integral is bounded by 
 $$A 
 \max(|1-t|^{-1+\lambda},|1-t|^{-1+\delta})\max(1,t^{-\lambda},t^{-p-\lambda}).$$
 \begin{lemma}
    Suppose $-1 <  p < 1, \,\, 0< \lambda, \delta <1, \mu< 0, p + \mu + 
    \lambda < 0$ and $ t > 1.$ Then 
    \be
    \int_{0}^{\infty} s^p (1+s)^{\mu} |1-s|^{-1+ 
    \lambda}|1-ts|^{-1+\delta}\,\,ds \leq A \max (|1-t|^{-1+\lambda}, 
    |1-t|^{-1+\delta}) \max (1,t^{-p })
   \ee 
   where $A$ is some constant independent of $t$.
   \end{lemma}
   Proof: The proof of this is almost identical to the previous lemma, and 
   we leave the details to the reader.
   \begin{lemma}
   Suppose $|x|<1, \Re{c}>0, \Re{(c-b)}>0,$ and $\Re{(c-a-b)}<0.$ Then 
   the hypergeometric function $F(a,b,c,x)$ satisfies the estimate 
   $$|F(a,b,c,x)| \leq A|1-x|^{\Re{(c-a-b)}}$$ with $A$ independent of $x$.
   \end{lemma} %end page 10
   Proof: The hypergeometric function satisfies the identity 
   $F(a,b,c,x)=(1-x)^{c-a-b}F(c-a,c-b,c,x).$ Using Euler's integral 
   formula for $F$, we have 
   \be
   F(c-a,c-b,c,x) = \frac{\Gamma (c)}{\Gamma (b) 
   \Gamma(c-b)}\int_{0}^{1}t^{c-b-1}(1-t)^{b-1}(1-tx)^{-c+a}\,dx.
   \ee
   The last integral is bounded by $\int_{0}^{1}t 
   ^{\Re{(c-b-1)}}(1-t)^{\Re{(a+b-c-1)}} dx$ or $\frac{\Gamma (\Re{ 
   (c-b)})\Gamma (\Re{(a+ b-c)}) }
   {\Gamma (\Re{a})}.$
   
   We next find an integral expression for the trace of 
   $(B_{\alpha}(\sigma))^{n}.$ We proceed informally at first
    and later state things rigorously.  Using (\ref{eqn2.11}) we can write this trace 
   as $$\int_{0}^{\infty} \ldots \int_{0}^{\infty} 
   \int_{0}^{1} \ldots \int_{0}^{1} \prod_{i=1}^{n}s_{i}x_{i}\sigma 
   (x_{i}/\alpha )J_{\nu}(x_{i}s_{i})J_{\nu}(x_{i}s_{i+1}) \,ds_{1} 
   \dots ds_{n} dx_{1} \ldots dx_{n}$$ where $s_{1+n} = s_{1}.$ Let 
   $\hat{\sigma}$ be the Mellin transform of $\sigma$ where $c > 0$. Then the above 
   becomes 
  $$
   \frac{1}{(2\pi i)^{n}}\int_{c-i\infty}^{c+i\infty} \ldots 
   \int_{c-i \infty}^{c+i \infty}
   \int_{0}^{\infty} \ldots \int_{0}^{\infty} 
   \int_{0}^{1} \ldots \int_{0}^{1} \prod_{i=1}^{n}  \{s_{i}x_{i}^{1-z_{i}} 
    J_{\nu}(x_{i}s_{i})J_{\nu}(x_{i}s_{i+1}) 
    \hat{\sigma}(z_{i})\}$$
    $$ \times  \alpha^{z_{1}+\dots z_{n}} \,ds_{1} 
   \dots ds_{n} dx_{1} \ldots dx_{n} dz_{1} \ldots dz_{n}. $$ %end page 11
   Now use the formula 
   $$\int_{0}^{\infty}x^{-\lambda}J_{\nu}(ax)J_{\nu}(bx)dx$$
   $$= \frac{(ab)^{\nu}\Gamma(\nu+\frac{1-\lambda}{2})}
   {2^{\lambda}(a+b)^{2\nu -\lambda +1}\Gamma (1+\nu) 
   \Gamma (1/2+\frac{\lambda}{2})} 
	  F(\nu+\frac{1-\lambda}{2},\nu+\frac{1}{2}; 2\nu 
	  +1;\frac{4ab}{(a+b)^{2}}),$$
	  where $F(a,b;c;z)$ is the hypergeometric function $_{2}F_{1}$, $n$ 
	  times in the integral to get the expression $$
   \frac{1}{(2\pi i)^{n}}\int_{c-i\infty}^{c+i\infty} \ldots 
   \int_{c-i \infty}^{c+i \infty}
   \int_{0}^{1} \ldots \int_{0}^{1} \alpha^{z_{1}+ \ldots +z_{n}}$$
    $$ \times \prod_{i=1}^{n} \hat{\sigma}(z_{i}) \frac{s_{i}^{2\nu 
    +1}\Gamma(\nu +1-z_{i}/2) 
    F(\nu + 1-z_{i}/2,\nu+\frac{1}{2}; 2\nu 
    +1;\frac{4s_{i}s_{i+1}}{(s_{i}+s_{i+1})^{2}}) }
    {2^{z_{i}-1}\Gamma (1+\nu)\Gamma (z_{i}/2)(s_{i}+s_{i+1})^{2\nu 
    -z_{i}+2}}$$ 
    $$ \times  ds_{1} 
   \dots ds_{n}   dz_{1} \ldots dz_{n}. $$
   Next we make the change of variables
   \begin{eqnarray*}
   s_{1} & = & s_{1}'\\
   s_{2} & = & s_{2}'s_{1}'\\
    & \vdots & \\
   s_{n} & = & s_{n}' \dots s_{1}'
   \end{eqnarray*}
   and the integral becomes
  
   $$\frac{1}{(2\pi i)^{n}}\int_{c-i\infty}^{c+i\infty} \ldots 
   \int_{c-i \infty}^{c+i \infty}
   \int_{0}^{1} \int_{0}^{\frac{1}{s_{1}}} \ldots \int_{0}^{\frac{1}{s_{1}\ldots 
   s_{n-1}}} (\alpha /2)^{z_{1}+ \ldots +z_{n}}
    2^{n} \hat{\sigma}(z_{n})
  \frac{\Gamma (\nu +1-z_{n}/2)}{\Gamma (1+\nu)
   \Gamma (z_{n}/2)}$$ 
   $$\times s_{1}^{z_{1}+ \ldots +z_{n}-1}(1+s_{n}\ldots 
   s_{2})^{-2\nu +z_{n}-2}
      F(\nu + 1-z_{n}/2,\nu+\frac{1}{2}; 2\nu 
    +1;\frac{4s_{n}\ldots s_{2}}{(1 + s_{n}\ldots s_{2})^{2}})$$
    $$\times \prod_{i=1}^{n-1}\{ \frac{\hat{\sigma}(z_{i}) \Gamma(\nu +1-z_{i}/2)}
    { \Gamma (1+\nu)\Gamma (z_{i}/2)}s_{i+1}^{2\nu +1+z_{i+1}+\ldots 
    z_{n-1}}(1+s_{i+1})^{-2\nu + z_{i}-2}$$  
     $$ \times F(\nu + 1-z_{i}/2,\nu+\frac{1}{2}; 2\nu 
    +1;\frac{4 s_{i+1}}{(1+s_{i+1})^{2}})\}  
      ds_{n} 
   \dots ds_{1}   dz_{1} \ldots dz_{n}. $$ Write the inside integral 
   as 
   $$\int_{0}^{1} \int_{0}^{\frac{1}{s_{1}}} \dots \int_{0}^{\frac{1}{s_{1}\ldots 
    s_{n-1}}} \ldots ds_{n} \dots d_{1} - \int_{0}^{1} \int_{0}^{\infty} 
    \dots \int_{0}^{\infty} \dots ds_{n} \dots ds_{1}$$
    $$ +\int_{0}^{1} \int_{0}^{\infty} 
    \dots \int_{0}^{\infty} \dots ds_{n} \dots ds_{1}.$$
    The last integral in the above sum, inserted in the main integral, is 
    the same as $\tr (B_{\alpha}(\sigma^{n})).$ After reversing the order 
    of integration to $ds_{1} \dots ds_{n},$ the first two terms combine 
    to yield limits of integration
    $$-\int_{0}^{\infty} \int_{0}^{\infty} \dots \int_{0}^{\infty} 
    \int_{\min(1,\frac{1}{s_{2}}, \dots ,\frac{1}{s_{2}\ldots 
    s_{n}})}^{1}$$ and then the first integration can be done. The result is 
    that 
    $$\tr(B_{\alpha}(\sigma))^{n} = \tr B_{\alpha}(\sigma^{n}) + 
    C(\sigma)$$
    where $C(\sigma)$ is given by the expression 
    $$\frac{-1}{(2\pi i)^{n}}\int_{c-i\infty}^{c+i\infty} \ldots 
   \int_{c-i \infty}^{c+i \infty}
     (\alpha /2)^{z_{1}+ \ldots +z_{n}} \prod_{i=1}^{n}\frac{\hat{\sigma}(z_{i})
     2 \Gamma(\nu +1-z_{i}/2)}
    { \Gamma (1+\nu)\Gamma (z_{i}/2)}$$ 
   $$\times \int_{0}^{\infty} \dots \int_{0}^{\infty}
     \frac{1-(\min(1,\frac{1}{s_{2}}, \dots ,\frac{1}{s_{2}\ldots 
    s_{n}}))^{z_{1} + \dots + z_{n}}}{z_{1}+z_{2}+\dots + z_{n}}$$
 $$\times (1+s_{n}\ldots 
   s_{2})^{-2\nu +z_{n}-2}
      F(\nu + 1-z_{n}/2,\nu+\frac{1}{2}; 2\nu 
    +1;\frac{4s_{n}\ldots s_{2}}{(1 + s_{n}\ldots s_{2})^{2}})$$
    $$\times \{\prod_{i=2}^{n}
     s_{i}^{2\nu +1+z_{i}+\ldots 
    z_{n-1}}(1+s_{i})^{-2\nu + z_{i-1}-2}  
      \times F(\nu + 1-z_{i-1}/2,\nu+\frac{1}{2}; 2\nu 
    +1;\frac{4 s_{i}}{(1+s_{i})^{2}})\} $$ 
    $$ \times ds_{n} 
   \dots ds_{2}   dz_{1} \ldots dz_{n}. $$ 
   We next write this integral as
   $$ \frac{-1}{(2\pi i)^{n}}\int_{c-i\infty}^{c+i\infty} \ldots 
   \int_{c-i \infty}^{c+i \infty}G(z_{i})\int_{0}^{\infty} \dots 
   \int_{0}^{\infty}H(z_{i};s_{i})ds_{n} 
   \dots ds_{2}   dz_{1} \ldots dz_{n}.$$ %end page 12
   The idea from here on out is to 
   evaluate this integral asymptotically using complex analysis. This 
   will be done in several stages and by breaking the integral into 
   several parts. To begin we first consider the interior integration 
   $$ \int_{0}^{\infty} \dots \int_{0}^{\infty} H(z_{i};s_{i}) ds_{2} \dots
   ds_{n}.$$
    Consider this as an integral over $R_{1}\cup R_{2}$ where $R_{1}$ is 
    a union of disjoint sets, $R_{1} = \cup_{i=2}^{n} U_{i}$ such that on 
    $U_{i}, s_{i}$ is bounded away from $1$ and where $R_{2}$ is the 
    complement of $R_{1}.$ 
    \begin{lemma}
    Suppose that $-2\nu -1 < 0$. The integral of $H(z_{i}; s_{i})$   
    over $U_{i}$ is bounded and the 
    $z_{i}$ 
    variables can be changed in such a way so that the integrated 
    function is analytic in a particular $z$ variable to the left of the imaginary 
    axis.
     \end{lemma}
    Proof: For convenience let $i=2$ (although the proof is the same for 
    any $i$) and let $z_{1}+ \ldots + z_{n} = 
    z_{1}^{'}$ with the other variables remaining the same. Suppose that 
    $\Re z_{i} = c $ for $i=3, \ldots ,n$ and that $c > 0. $ Suppose also 
    that $\Re z_{1}^{'} = b$ with $|b| < c.$ We now refer to $z_{1}'$ as $z$.
    Our goal is to show that this integral is bounded 
    and that as a function of $z$ is analytic to the left of the 
    imaginary axis. By repeated application of Lemma 8, we can say that the integral is 
    bounded by a constant times 
    $$\int_{|s_{2}-1| \geq B}\int_{0}^{\infty} \ldots \int_{0}^{\infty}
     |\frac{1-(\min(1,\frac{1}{s_{2}}, \dots ,\frac{1}{s_{2}\ldots 
    s_{n}}))^{z_{1} + \dots + z_{n}}}{z_{1}+z_{2}+\dots + z_{n}}|$$
    $$\times \prod_{i=3}^{n} (1 + s_{i})^{-2\nu - 1}s_{i}^{
     (n-i)c} |1 -s_{i}|^{c -1}$$
    $$\times s_{2}^{b-2c}(1+s_{2})^{-2\nu-1}|1-s_{2}|^{b-1}
    |1-s_{2} \ldots s_{n}|^{c-1}$$
    $$\times ds_{2} \ldots ds_{n}.$$ This is valid as long as
     $ 2\nu +1 > 0 $ and $\Re z_{i} - 1/2 < 0$, 
    which is the case here if we assume that $c$ is small enough.
    Next, we estimate 
    $$|\frac{1-(\min(1,\frac{1}{s_{2}}, \dots ,\frac{1}{s_{2}\ldots 
    s_{n}}))^{z_{1} + \dots + z_{n}}}{z_{1}+z_{2}+\dots + z_{n}}|$$ by 
    using the fact that $|1- x^{z}| \leq \max |z|x^{\Re z'}|\ln x|$ where 
    the max is taken over the $z'$ values on a line connecting $0$ and $z$, $x$ 
    between $0$ and $1$. Thus, $|1- x^{z}| \leq K|z|x^{\Re z} 
    x^{-\epsilon}$ for some positive $\epsilon$ chosen shortly. Inserting this in the 
    integral we have that the integral is bounded by a constant times
     $$\int_{|s_{2}-1| \geq B}\int_{0}^{\infty} \ldots \int_{0}^{\infty}
      \sum_{j=2}^{n}\{\max(1,(s_{2}\ldots s_{j})^{\epsilon} 
      (s_{2}\ldots s_{j})^{-b+\epsilon})\}$$
    $$\times \prod_{i=3}^{n} (1 + s_{i})^{-2\nu - 1}s_{i}^{
     (n-i)c} |1 -s_{i}|^{c -1}$$
    $$\times s_{2}^{b-2c}(1+s_{2})^{-2\nu-1}|1-s_{2}|^{b-1}
    |1-s_{2} \ldots s_{n}|^{c-1}$$
     $$\times ds_{2} \ldots ds_{n}.$$ %end page 13
     The reason for both terms in the 
    ``max" part of the integral is that $b$ could be either positive or 
    negative. Now lets begin with the $s_{n}$ integration.  
  Then the first interior
    integral has the 
    form $$\int _{0}^{\infty} s_{n}^{p}(1+s_{n})^{-2\nu-1}|1-s_{n}|^{c-1}
    |1-s_{2}\ldots s_{n}|^{-1 +c}ds_{n}.$$ The value for $p$ is 
    either $\pm a $ where $a = |b -\epsilon| < c$. The next step 
    is to apply Lemmas 6 and 7. We use $\lambda = \delta = c$ and $p$ as 
    above. The result is that this integral is bounded by a constant times
    $$ |1-s_{2}\ldots s_{n-1}|^{c-1} \times \max (1, (s_{2}\ldots 
    s_{n-1})^{-c}, (s_{2} \ldots s_{n-1})^{-c-p}(s_{2} \ldots 
    s_{n-1})^{-p}).$$ We collect powers and use the lemmas twice  
    with respect to the $s_{n-1}$ integration 
    and powers of $p = \pm(2c)$. At the next 
    integration step the powers of $p = \pm 3c$ and 
    so on until we arrive at the $s_{2}$ integration. Here we will have $$
    \int_{|s_{2}-1|\geq B}s_{2}^{p}|1-s_{2}|^{q}|1+s_{2}|^{-2\nu-1}ds_{2}$$
    where $p$ and $q$ are appropriate powers.  
    These integrals satisfy all the conditions necessary for the lemmas 
    as long as $c$ and $b$ are small enough. We will have at most $2^{n}$ 
    integrals in this process. Hence the integral of $H$ over $U_{i}$ is 
    analytic in the z variable in a strip $|\Re z|<c$ by the application of Morera's Theorem and
    Fubini's Theorem.
    
    We remark here that this proof also is easily modified to show that the interchange of 
    integrals done at the beginning of the section are valid and the 
    expression $C(\sigma)$ is the one of interest.

   \begin{lemma}
   Suppose that $\sigma$ has $[\nu]+2$ derivatives all in $L_{1}$ 
     and that that $-2\nu - 1 < 0$. Then the integral 
   $$\int_{c-i\infty}^{c+i\infty} \ldots \int_{c-i\infty}^{c+i\infty} G(z_{i})
   \int \ldots \int_{R_{1}} H(z_{i};s_{i})ds_{2} 
   \ldots ds_{n}dz_{1} \ldots dz_{n}$$  
  is $O(\alpha^{-\delta})$ where $\delta > 0.$  
  \end{lemma}
  Proof: 
  Note that the condition in the hypothesis implies that
      \be
   \int_{c-i\infty}^{c+\infty}|\hat{\sigma}(z)||z|^{\nu+1/2} < \infty.
    \label{eqn4.1}
   \ee
  
  We first replace the inside integral with a sum of integrals 
  over $U_{i}.$ For each of these we change variables as in the last 
  lemma. We can then perform the integration over the $z$ variable by 
  moving it to a line to the left of the imaginary axis. Thus we have that 
  each of these integrals is bounded by a constant times 
  $$  \frac{\alpha^{b}}{(\pi )^{n}2^{b}}\int_{b-i\infty}^{b+i\infty}\int_{c-i\infty}^{c+i\infty} \ldots 
   \int_{c-i \infty}^{c+i \infty}
      | \prod_{i=2}^{n}\frac{\hat{\sigma}(z_{i})
      \Gamma(\nu +1-z_{i}/2)}
    { \Gamma (1+\nu)\Gamma (z_{i}/2)}|$$ 
   $$\times |\frac{\hat{\sigma}(z-\sum_{z_{j} \neq z}z_{j})\Gamma(\nu + 1-(z-\sum 
   _{z_{j}\neq z}z_{j})/2)}{\Gamma(1+\nu) \Gamma((z -\sum_{z_{j}\neq 
   z}z_{j})/2)}|dzdz_{2} \ldots 
   dz_{n}.$$ %end page 14
   This last integral is bounded by a product of integrals all of the 
   form $$\int_{c-i\infty}^{c+i\infty}|\hat{\sigma}(z)||\frac{\Gamma (\nu +1 
   -z/2)}{\Gamma (z/2)}|dz$$ and these in turn are bounded by 
   (\ref{eqn4.1}) using the basic asymptotics properties of the 
   Gamma function.
    
    We now turn our attention to the region $R_{2}.$ To begin we make 
    another 
    change of variables.
    \begin{eqnarray}
    \frac{1}{s_{2}} &= &1-s_{2}'\\
    \frac{1}{s_{2}s_{3}}& = &1-s_{2}'-s_{3}'\\
     & \vdots &\\
    \frac{1}{s_{2}\ldots s_{n}}& = & 1-s_{2}'-s_{3}'-\ldots -s_{n}'.
    \end{eqnarray}  
    Under the change of variables, the region $R_{2}$ is transformed to a 
    region $R_{3}$ which can be assumed to be a symmetric region 
    containing the origin, and where the sum $$|s_{2}+ \ldots +s_{j}| \leq 
    a <1$$ (we drop the ``primes" again) for some $a.$ Notice that the exact form of $R_{1}$ was 
    unnecessary in the previous computation. Thus the integral over 
    $R_{2}$ is transformed to $$ \raisebox{-4mm}{$\displaystyle 
    \int \ldots \int \atop \hspace{-2mm} 
    R_{3} $}I(z_{i};s_{i})ds_{2}\ldots ds_{n}$$ where
    $$I(z_{i};s_{i}) = \frac{1-(1-\max(0,s_{2},\ldots ,s_{2}+\ldots 
    +s_{n}))^{z_{1}+\ldots +z_{n}}}{z_{1}+ \ldots 
    +z_{n}}$$
    $$\times|s_{2}|^{z_{2}-1}\ldots |s_{n}|^{z_{n-1}-1}|s_{2}+\ldots 
    +s_{n}|^{z_{n}-1}$$
    $$\times f(s_{2},\ldots ,s_{n},z_{1},\ldots , z_{n})$$
    where the function $f$ is smooth in the $s$ variables.
    
    The following lemmas will help keep track of the contribution of the 
    $R_{3}$ integral.
    \begin{lemma}
    Suppose $\Re z_{i} = c, 0<c<1,$ for $i\geq 3.$ Then the integral 
    $$ \raisebox{-4mm}{$\displaystyle 
    \int \ldots \int \atop \hspace{-2mm} 
    R_{3} $} |s_{2}|^{z_{1}+1}|s_{3}|^{z_{2}-1} \ldots 
    |s_{n}|^{z_{n-1}-1}|s_{2}+\ldots s_{n}|^{z_{n}-1}ds_{2} \ldots ds_{n}$$
    can be thought of as an analytic function in the $z_{1}$ variable that 
    can be extended to a strip containing the imaginary axis.
    \end{lemma}
    Proof: First note that the following integral with $z$ and $w$ real and between 
    zero and one satisfies
    $$\int_{a}^{b}|x|^{z-1}|x+y|^{w-1}dx \leq A|y|^{z+w-1}$$ where the 
    constant only depends on the $z$ and $w$ variable. A repeated 
    application of this estimate in the above integral yields a final 
    integration of $$ \int_{a}^{b}|s_{2}|^{\Re z_{1}+(n-2)c}ds_{2}. $$
    Thus, once again the analytic continuation 
    argument holds. %end page 15
    \begin{lemma}
    Suppose $\Re z_{i} = c, 0<c<1,$ for $ i\geq 2$ and $\Re z_{1}=d$. Then the integral
    $$ \raisebox{-4mm}{$\displaystyle 
    \int \ldots \int \atop \hspace{-2mm} 
    R_{3} $}
    |s_{2}|^{z_{1}}|s_{3}|^{z_{2}}|s_{4}|^{z_{3}-1} \ldots 
    |s_{n}|^{z_{n-1}-1}|s_{2}+\ldots s_{n}|^{z_{n}-1}ds_{2} \ldots ds_{n}$$
    can be thought of as an analytic function in the $z_{1}$ variable 
    that can be extended to a strip containing the imaginary axis.
    \end{lemma}
    Proof: We begin the integration just as in the previous integral. 
    After $n-3$ integrations we arrive at an integral with an estimate 
    of the form
    $$\int_{a}^{b}\int_{a}^{b}|s_{2}|^{d}|s_{3}|^{c}|s_{3}+s_{2}|^{ 
    (n-3)c-1}ds_{2}ds_{3}.$$
   We can estimate this by looking at three integrals
   $$\int_{a}^{b}\int_{-1}^{1}|s_{2}|^{ d+(n-2)c}|s_{3}|^{c}|s_{3}+1|^{ 
    (n-3)c-1}ds_{3}ds_{2},$$
    $$\int_{a}^{b}\int_{1}^{b/s_{2}}|s_{2}|^{ d+(n-2)c}|s_{3}|^{c}|s_{3}+1 |^{ 
    (n-3)c-1}ds_{3}ds_{2},$$ and
    $$\int_{a}^{b}\int_{a/s_{2}}^{-1}|s_{2}|^{ d+(n-2)c}|s_{3}|^{c}|s_{3}+1 |^{ 
    (n-3)c-1}ds_{3}ds_{2}.$$
    We can say, for example, that the last integral is less than a 
    constant times $$\int_{a}^{b}|s_{2}|^{d-(n-4)c}ds_{2}$$ and thus is 
    finite for $\Re z_{1}$ in a strip about the imaginary axis. The other 
    two integrals are handled in the same manner. So by our standard 
    argument the analytic extension is defined.
    
    Now let us return to our function $I(z_{i};s_{i})$. We can write the 
    expression $$\frac {1-(1-\max (0,s_{2},\ldots, s_{2}+\ldots 
    +s_{n}))^{z_{1}+\ldots z_{n}}}{z_{1}+\ldots + z_{n}}$$ as 
    $$\max (0,s_{2},\ldots, s_{2}+\ldots s_{n}) + (\max (0,s_{2},\ldots, 
    s_{2}+\ldots s_{n}))^{2}
    \times g(z_{1}+\ldots + z_{n}, s_{2}, \ldots s_{n})$$ where the last 
    function is a bounded continuous function in the variables.
    
    \begin{lemma} The contribution of $$\frac{-1}{(2\pi i)^{n}}
    \int_{c-i\infty}^{c+i\infty} \ldots 
   \int_{c-i \infty}^{c+i \infty}G(z_{i})\int  \dots 
   \int_{ R_{2}} (\max (0,s_{2},\ldots, s_{2}+\ldots s_{n}))^{2}$$
    $$\times|s_{2}|^{z_{1}-1}\ldots |s_{n}|^{z_{n-1}-1}|s_{2}+\ldots 
    +s_{n}|^{z_{n}-1}$$
    $$\times f(s_{2},\ldots, s_{n},z_{1}, \ldots , z_{n})
    g(z_{1}+\ldots + z_{n}, s_{2}, \ldots s_{n})$$
    $$\times ds_{n} 
   \dots ds_{2}   dz_{1} \ldots dz_{n}$$ is $O(\alpha^{-\delta}).$ 
   \end{lemma}
   Proof: We simply consider the set where say $s_{1} + s_{2} +\ldots 
   + s_{j}$ is the maximum of the the terms. We then expand the square so 
   that we have a term of the form $s_{i}s_{k}$. We then apply the above 
   lemmas after an appropriate re-ordering of the variables and the lemma 
   holds.
   
   The next step is to replace the function $f$ in the expression for 
   $I(z_{i};s_{i})$ with the first term of its Taylor expansion. This 
   expansion gives an ``extra" $s_{i}$ (combined with the ones from the 
   $\max(0, s_{2}, \ldots ,s_{2}+\ldots +s_{n})$) term in the estimates 
   which, as the 
   above lemmas show, is all we need to show that this part of the integral 
   does not contribute in the asymptotic expansion.
   
   So we are finally at the one critical term that gives a contribution 
   in the expansion. This term is 
   $$
   \frac{-1}{(2\pi i)^{n}}\int_{c-i\infty}^{c+i\infty}\ldots \int_{c-i\infty}^{c+i\infty}
    G(z_{i})
   \int \ldots \int_{R_{3}}\max(0, s_{2}, \ldots ,s_{2}+\ldots +s_{n})$$
   $$\times |s_{2}|^{z_{1}-1} \ldots |s_{n}|^{z_{n-1}-1}|s_{2}+\ldots 
   +s_{n}|^{z_{n}-1}$$
   $$\times f(0,0,\ldots,0,z_{1},\ldots ,z_{n})ds_{n}\ldots ds_{2}dz_{1} 
   \ldots d_{n}.$$
   We can easily compute $f(0,0,\ldots,0,z_{1},\ldots, z_{n})$ to see 
   that it equals $$\prod_{1}^{n}\frac{2^{-2\nu -1}\Gamma (2\nu 
   +1)\Gamma (-z_{i}/2+1/2)}{\Gamma (\nu +1/2)\Gamma (\nu +1-z_{i}/2)}.$$
  We can simplify further using the formula for $G(z_{i})$ and the 
  duplication formula for the Gamma function to arrive at
  $$
  C(\sigma)   = \frac{-1}{(2\pi i)^{n}}\int_{c-i\infty}^{c+i\infty}\ldots \int_{c-i\infty}^{c+i\infty}
    (\alpha /2)^{z_{1}+\ldots +z_{n}} \pi^{-n/2} 
    \prod_{1}^{n}\frac{\hat{\sigma} (z_{i}) \Gamma (-z_{i}/2+1/2)}{\Gamma 
    (z_{i}/2)}$$
$$ \times\int \ldots \int_{R_{3}}\max(0, s_{2}, \ldots ,s_{2}+\ldots 
   +s_{n})  $$
 \be
   \times |s_{2}|^{z_{1}-1} \ldots |s_{n}|^{z_{n-1}-1}|s_{2}+\ldots 
   +s_{n}|^{z_{n}-1}
    ds_{2} \dots ds_{n}dz_{1} \ldots dz_{n} + O(\alpha^{-\delta})
   . \label{eqn4.2}
   \ee
 
 Notice that this expression is now independent of $\nu.$ Our final steps 
 are to compute the contribution from the above integral and we, by the way, finally 
 have an integral which will yield a contribution. We begin with a 
 well-known 
 identity due to Mark Kac, which was used originally to prove the 
 continuous analogue of the Strong Szeg\"{o} Limit Theorem. It reads
 $$\sum_{\sigma}\max(0, a_{\sigma_{1}}, a_{\sigma_{1}}+a_{\sigma_{2}}, 
 \ldots ,a_{\sigma_{1}}+ \ldots +a_{\sigma_{n}}) = 
 \sum_{\sigma}\sum_{k=1}^{n}a_{\sigma_{1}} \theta ( a_{\sigma_{1}}+\ldots 
 +a_{\sigma_{k}})$$ where $\theta (x) = 1$ if $x > 0$ and $\theta (x) = 
 0$ otherwise and the sums are taken over all permutations in $n$ variables.
 
 Because of this identity we can rewrite the integral in (\ref{eqn4.2}) as 
 $$    \sum_{j=2}^{n} \frac{-1}{(2\pi i)^{n}}\int_{c-i\infty}^{c+i\infty}\ldots \int_{c-i\infty}^{c+i\infty}
    (\alpha /2)^{z_{1}+\ldots +z_{n}} \pi^{-n/2} 
    \prod_{1}^{n}\frac{\hat{\sigma} (z_{i}) \Gamma (-z_{i}/2+1/2)}{\Gamma 
    (z_{i}/2)}$$
 \be
 \times\int \ldots \int_{R_{3}\cap \{s_{2}+\ldots + s_{j}>0\}
} s_{2} 
      |s_{2}|^{z_{1}-1} \ldots |s_{n}|^{z_{n-1}-1}|s_{2}+\ldots 
   +s_{n}|^{z_{n}-1}
    ds_{2} \dots ds_{n}dz_{1}. \ldots dz_{n}. \label{eqn4.3}
    \ee
    It is straightforward to see how this identity can be used if the 
    integrand is symmetric in the variables. In our case, the integrand 
    is not obviously symmetric in the variables, but can always be made 
    so by changing the $z$ variables. Thus we can apply the identity.
    %end page 17
    
    We once again consider the inner integral and call $z = z_{1}+\ldots 
    +z_{n}$ leaving the other variables as is, and show how this inner 
    integral can be thought of as analytic in $z$ in a strip containing the 
    imaginary axis. The difference is that in this case there will be a
    pole at $z = 0.$
 
 Now we suppose that $j >2$. For $j=2$ the following computation is 
 almost  identical and the conclusion is the same.
  Let us rewrite the inner integral in (\ref{eqn4.3}) as two integrals
 $$\int_{0}^{b} \raisebox{-4mm}{$\displaystyle 
    \int \ldots \int \atop \hspace{-2mm} 
   B $} s_{2} 
      |s_{2}|^{z -z_{2}-z_{3}-\ldots z_{n}-1} \ldots |s_{n}|^{z_{n-1}-1}|s_{2}+\ldots 
   +s_{n}|^{z_{n}-1}
    ds_{n} \dots ds_{2} $$
    $$ + \int_{-b}^{0}\raisebox{-4mm}{$\displaystyle 
    \int \ldots \int \atop \hspace{-2mm} 
   B $} s_{2} 
      |s_{2}|^{z -z_{2}-z_{3}-\ldots z_{n}-1} \ldots |s_{n}|^{z_{n-1}-1}|s_{2}+\ldots 
   +s_{n}|^{z_{n}-1}
    ds_{n} \dots ds_{2} $$ where $B$ is some $n-2$ dimensional set.
    In the first (the computations for the second integral being almost 
    identical) of these we make yet another change of variables:
    \begin{eqnarray*}
  s_{3} & = & s_{3}'s_{2}'\\
    & \vdots & \\
   s_{n} & = & s_{n}' s_{2}'
   \end{eqnarray*}
  
   to arrive at $$\int_{0}^{b} s_{2}^{z -1}
   \raisebox{-4mm}{$\displaystyle 
    \int \ldots \int \atop \hspace{-2mm} 
   B/s_{2} $}  
      |s_{3}|^{z_{2}-1}\ldots   |s_{n}|^{z_{n-1}-1}|1+s_{3}+\ldots 
   +s_{n}|^{z_{n}-1}
    ds_{n} \dots ds_{3} ds_{2}. $$
    The original set $R_{3}$ was chosen to be symmetric and contain 
    the origin. So here we chose it to be something convenient, say a
    cube $C$ with size length $l$. With 
    this choice we can write $B/s_{2}$ as $C/s_{2}\cap \{s_{3}+\ldots + 
    s_{n} +1>0\}.$ Next integrate by parts with respect to the $s_{2}$ 
    variable. The result is that the above integral becomes:
    
    $$s_{2}^{z}k(s_{2})
    - \int_{0}^{b}s_{2}^{z}  d/ds_{2}(k(s_{2}))ds_{2}$$ where $$k(s_{2})
    =\raisebox{-4mm}{$\displaystyle 
    \int \ldots \int \atop \hspace{-2mm} 
   B/s{2} $}  
      |s_{3}|^{z_{2}-1} \ldots   |s_{n}|^{z_{n-1}-1}|1+s_{3}+\ldots 
   +s_{n}|^{z_{n}-1}
    ds_{n} \dots ds_{3}.$$ The function $k(s_{2})$ has a derivative given 
    by the formula $$ k'(s_{2}) = -s_{2}^{-1}\int_{D}f(s_{3},\ldots, 
    s_{n})\ (\vec{n} \cdot s_{2}^{-1}(s_{3},\ldots, s_{n}))dS$$ where $D$ 
    is the boundary of the set $C/s_{2}$ which lies in the half-space 
    defined by $\{s_{3}+\ldots + s_{n} +1>0\},$ the vector $ \vec{n}$ is 
    the outward normal to the surface, the function $f$ is simply the one 
    given in the above integral restricted to the surface, and $dS$ is 
    surface measure. We can estimate the derivative of $k(s_{2})$ on any boundary edge 
    to be at most a constant times $s_{2}^{(n-2)c}$ for $\Re z_{i} = c$. 
    Thus we have proved the following: %end page 18
    \begin{lemma}
    The function of $z$ defined by
    $$\int_{0}^{b}\raisebox{-4mm}{$\displaystyle 
    \int \ldots \int \atop \hspace{-2mm} 
   B $} s_{2} 
      |s_{2}|^{z -z_{2}-z_{3}-\ldots z_{n}-1} \ldots |s_{n}|^{z_{n-1}-1}|s_{2}+\ldots 
   +s_{n}|^{z_{n}-1}
    ds_{n} \dots ds_{2} $$
    $$ + \int_{-b}^{0}\raisebox{-4mm}{$\displaystyle 
    \int \ldots \int \atop \hspace{-2mm} 
   B $} s_{2} 
      |s_{2}|^{z -z_{2}-z_{3}-\ldots z_{n}-1} \ldots |s_{n}|^{z_{n-1}-1}|s_{2}+\ldots 
   +s_{n}|^{z_{n}-1}
    ds_{n} \dots ds_{2} $$ is analytic in a strip containing the 
    imaginary axis except at the point $z =0$. Further, the contribution 
    of this integral with the $z$ integration moved to a line to the left 
    of the axis is given by the residue at $z = 0$ plus 
    $O(\alpha^{-\delta}).$
    \end{lemma}
    
    We note here that there are no other poles given our conditions on 
    $\sigma$ , ( \ref{eqn4.1}) and the formula for $G(z_{i})$.
    
      For $j>2,$  the above computation also shows exactly what the residue is, namely:
    $$ \raisebox{-4mm}{$\displaystyle 
    \int \ldots \int \atop \hspace{-2mm} 
   \R^{n-2} \cap \{ s_{3}+ \ldots + s_{j} >-1\} $}   
      |s_{3}|^{z_{2}-1} \ldots   |s_{n}|^{z_{n-1}-1}|1+s_{3}+\ldots 
   +s_{n}|^{z_{n}-1}
    ds_{n} \dots ds_{3}$$
    $$ - \raisebox{-4mm}{$\displaystyle 
    \int \ldots \int \atop \hspace{-2mm} 
   \R^{n-2} \cap \{ s_{3}+ \ldots + s_{j} >-1\} $}    
      |s_{3}|^{z_{2}-1} \ldots   |s_{n}|^{z_{n-1}-1}|-1+s_{3}+\ldots 
   +s_{n}|^{z_{n}-1}
    ds_{n} \dots ds_{3}.$$
    
    To find an explicit formula for this integral we start with the 
    following formula that can be easily proved using formulas for the Beta 
    function.
    
    For $$0 < \Re p, \Re q < 1, \Re (p+q) <1 $$
    \be 
  \int_{-\infty}^{\infty}|x|^{p-1}|x+y|^{q-1}dx = 
    |y|^{p+q -1}\frac{2\Gamma (p) \Gamma (q) \cos(\pi p/2) \cos (\pi q/2)}{\Gamma (p+q) 
    \cos ((p+q)\pi/2)}.
    \ee
    Define $t(p,q) $ to be $$\frac{2\Gamma (p) \Gamma (q) \cos(\pi p/2) \cos (\pi q/2)}{\Gamma (p+q) 
    \cos ((p+q)\pi /2)}.$$ The residue is then ($B$ is the Beta function)
    $$ B(z_{2}+\ldots + z_{j-1},z_{j}+\ldots + z_{n}) \prod_{k=j}^{n-1} t(z_{k},z_{n}+ \ldots +z_{k+1}) 
    \prod_{k=2}^{j-2}t(z_{k},z_{k+1}+ \ldots  +z_{j-1}).$$ We leave 
    this as an exercise to the reader.
    For $j=2$ the residue can also be easily computed using the 
    definition of $t(p,q)$ and it is seen to be
    $$\prod_{k=2}^{n-1}t(z_{k},z_{n}+\ldots +z_{k+1}).$$
    
    Combining all of the above results we are left with the following 
    theorem.
    \begin{theorem}
    Suppose $\sigma$ has $[\nu]+2$ continuous derivatives in $L_{1}$. Then 
    $$\tr(B_{\alpha}(\sigma))^{n} = \tr B_{\alpha}(\sigma^{n}) + 
    C(\sigma)$$ where 
    $$ C(\sigma) = \frac{-1}{\pi^{2}}\sum_{j=1}^{n-1}\frac{1}{j} 
    \int_{0}^{\infty}x C(\sigma^{j})(x) 
    C(\sigma^{n-j})(x)dx + o(1).$$ 
    \end{theorem}
    Proof: Recall we were computing the integral %end page 19
     $$    \sum_{j=2}^{n} \frac{-1}{(2\pi i)^{n}}\int_{c-i\infty}^{c+i\infty}\ldots \int_{c-i\infty}^{c+i\infty}
    (\alpha /2)^{z_{1}+\ldots +z_{n}}  \pi^{-n/2} 
    \prod_{1}^{n}\frac{\hat{\sigma} (z_{i}) \Gamma (z_{i}/2+1/2)}{\Gamma 
    (z_{i}/2)}$$
 \be
 \times\int \ldots \int_{R_{3}\cap \{s_{2}+\ldots + s_{j}>0\}
} s_{2} 
      |s_{2}|^{z_{1}-1} \ldots |s_{n}|^{z_{n-1}-1}|s_{2}+\ldots 
   +s_{n}|^{z_{n}-1}
    ds_{2} \dots ds_{n}dz_{1}. \ldots dz_{n}. \label{eqn4.4}
    \ee 
    For each $j$ we rename the variables and compute the residue  as 
    above. For $j>2$ 
  the residue is  
     $$      \frac{-1}{(2\pi i)^{n-1}}\int_{c-i\infty}^{c+i\infty}\ldots 
     \int_{c-i\infty}^{c+i\infty}
     \pi^{-n/2} 
    \prod_{2}^{n}\frac{\hat{\sigma} (z_{i}) \Gamma (-z_{i}/2+1/2)}{\Gamma 
    (z_{i}/2)}$$
    $$\times \frac{\hat{\sigma} (-z_{2}- \ldots -z_{n}) \Gamma ((z_{2}+ \ldots 
    +z_{n})/2+1/2)}{\Gamma 
    ( (-z_{2}- \ldots -z_{n})/2)}
 B(z_{2}+\ldots + z_{j-1},z_{j}+\ldots + z_{n})$$
    \be
    \times \prod_{k=j}^{n-1} t(z_{k},z_{n}+ \ldots +z_{k+1}) 
    \prod_{k=2}^{j-2}t(z_{k},z_{k+1}+ \ldots  +z_{j-1})dz_{2}\ldots dz_{n}.
    \ee \label{eqn4.5}
    Notice that 
    $$ t(p,q)t(p+q,r)=2^{2} \frac{\Gamma (p) \Gamma(q) \Gamma(r) \cos (p) 
    \cos (q) \cos (r)}{\Gamma (p+q+r) \cos ((p+q+r)\pi /2)}.$$
      Using this identity in (\ref{eqn4.5}) we have that the above 
      integral is 
      $$     \frac{-1}{(2\pi i)^{n-1}}\int_{c-i\infty}^{c+i\infty}\ldots 
     \int_{c-i\infty}^{c+i\infty}
     2^{n-3}\pi^{-n/2} 
    \prod_{i=2}^{n}\frac{\hat{\sigma} (z_{i}) \Gamma (z_{i}) \cos(z_{i}\pi 
    /2) \Gamma (-z_{i}/2+1/2)}{\Gamma 
    (z_{i}/2)}$$
    \be
    \times \frac{\hat{\sigma} (-z_{2}- \ldots -z_{n}) \Gamma ((z_{2}+ \ldots 
    +z_{n})/2+1/2) }{\Gamma 
    ( (-z_{2}- \ldots -z_{n})/2) \Gamma (z_{2}+\ldots +z_{n}) \cos 
    ((z_{2}+\ldots +z_{j-1})\pi /2) \cos ((z_{j}+\ldots +z_{n})\pi /2)}dz_{2}\ldots dz_{n}.
   \label{eqn4.6}
   \ee
   From the duplication formula for the Gamma function, this can be 
   simplified to  
   $$     \frac{-1}{(2\pi i)^{n-1}}\int_{c-i\infty}^{c+i\infty}\ldots 
     \int_{c-i\infty}^{c+i\infty}
     2^{-2}\pi^{-1} 
    \prod_{i=2}^{n} \hat{\sigma} (z_{i}) $$
    \be
    \times \frac{\hat{\sigma} (-z_{2}- \ldots -z_{n}) (z_{2}+\ldots +z_{n})  \sin ((z_{2}+\ldots 
    +z_{n})\pi /2)}{  \cos 
    ((z_{2}+\ldots +z_{j-1})\pi /2) \cos ((z_{j}+\ldots +z_{n})\pi /2)}dz_{2}\ldots dz_{n}.
   \label{eqn4.7}
   \ee
   Now we change variables with $$z_{j-1}=z_{2}+\ldots 
   +z_{j-1},z_{n}=z_{j}+\ldots +z_{n}$$   and the above integral becomes
    $$     \frac{-1}{(2\pi i)^{n-1}}\int_{c-i\infty}^{c+i\infty}\ldots 
     \int_{c-i\infty}^{c+i\infty}
     2^{-2}\pi^{-1} 
    (\prod_{i=2}^{j-2} \hat{\sigma} (z_{i}) )
    \hat{\sigma} (z_{j-1}-\ldots -z_{2}) \prod_{i=j}^{n-1} 
    \hat{\sigma} (z_{i})$$
    \be
    \times \hat{\sigma}(z_{n}-\ldots -z_{j}) 
    \hat{\sigma}(-z_{j-1}-z_{n})
    (z_{j-1} +z_{n}) \frac{ \sin ((z_{j-1} 
    +z_{n})\pi /2)}{  \cos 
    ( z_{j-1}\pi /2) \cos ( z_{n}\pi /2)}dz_{2}\ldots dz_{n}.
 \label{eqn4.8}
   \ee
   The convolution theorem for the Mellin transform shows that this can 
   be reduced to the integral 
     $$   \frac{-1}{(2\pi i)^{2}}\int_{c-i\infty}^{c+i\infty} 
     \int_{c-i\infty}^{c+i\infty}
     2^{-2}\pi^{-1}  
    \hat{\sigma^{j-2}} (z_{j-1})   
     \hat{\sigma^{n-j}}(z_{n}) $$
     \be
    \times \hat{\sigma}(-z_{j-1}-z_{n})
    (z_{j-1} +z_{n}) \frac{ \sin ((z_{j-1} 
    +z_{n})\pi /2)}{  \cos 
    ( z_{j-1}\pi /2) \cos ( z_{n}\pi /2)}dz_{j-1} dz_{n}.
   \label{eqn4.9}
  \ee 
   
  Notice this can also be written as
   $$    \frac{-1}{(2\pi  i)^{2}}\int_{c-i\infty}^{c+i\infty} 
     \int_{c-i\infty}^{c+i\infty}
      2^{-2}\pi^{-1} \hat{\sigma^{j-2}} (z_{j-1})   
     \hat{\sigma^{n-j+1}}(z_{n}) $$
     \be
    \times \hat{\sigma}(-z_{j-1}-z_{n})
    (z_{j-1} +z_{n}) \frac{ \sin (z_{j-1} 
    \pi /2)}{  \cos 
    ( z_{j-1}\pi /2)  }dz_{j-1} dz_{n}  
   \ee
    $$   - \frac{1}{(2\pi  i)^{2}}\int_{c-i\infty}^{c+i\infty} 
     \int_{c-i\infty}^{c+i\infty} 2^{-2}\pi^{-1} 
     \hat{\sigma^{j-2}} (z_{j-1})   
     \hat{\sigma^{n-j+1}}(z_{n}) $$
     \be
    \times \hat{\sigma}(-z_{j-1}-z_{n})
    (z_{j-1} +z_{n}) \frac{ \sin (z_{ n} 
    \pi /2)}{  \cos 
    ( z_{n}\pi /2)  }dz_{j-1} dz_{n}. \label{eqn4.10} 
   \ee

   Before we proceed further we need three formulas from the theory of 
   Mellin transforms. These are
   $$\mbox{the Mellin transform of} \int_{x}^{\infty} \phi (x)dx = 
   z^{-1}\Phi (z+1) $$ where $\Phi$ is the transform of $\phi$,
    $$ \mbox{the Mellin transform of } \,\,\,\, x \phi'(x) = 
   -z\Phi (z) $$ where $\Phi$ is the transform of $\phi$, and finally 
   $$\frac{2}{\pi}\int_{0}^{\infty} x C(\phi)(x) C(\psi)(x)dx = \frac{1}{2\pi 
   i}\int_{c-i\infty}^{c+\infty}\Phi(z) \Psi(-z) z \tan (z\pi /2)dz.$$
   These can be found in any standard table of transforms, although the 
   third requires a straightforward computation combined with the 
   convolution theorem.
   
   So now we apply the second formula along with convolution with respect 
   to the $ z_{n}$ variable and we have for each $2<j<n$  
   \be
       \frac{1}{8\pi^{2}i }\int_{c-i\infty}^{c+i\infty} 
    \widehat{\sigma ^{j-2}} (z_{j-1})   
     \widehat{x \sigma^{n-j+1} \sigma ' }(-z_{j-1}) 
    \frac{ \sin (z_{j-1} 
    \pi /2)}{  \cos 
    ( z_{j-1}\pi /2)  }dz_{j-1}  
  \ee
  \be
    +\frac{1}{8\pi^{2} i }\int_{c-i\infty}^{c+i\infty} 
    \widehat{\sigma ^{n-j+1}} (z_{j-1})   
     \widehat{x \sigma^{ j-2} \sigma ' }(-z_{j-1}) 
    \frac{ \sin (z_{j-1} 
    \pi /2)}{  \cos 
    ( z_{j-1}\pi /2)  }dz_{j-1}  .
   \label{eqn4.11}
  \ee
  Next apply the first formula after inserting a factor of 
  $z_{j-1}/z_{j-1}$ to write the above as
   \be
      \frac{1}{2 \pi^{2}} \int_{0}^{\infty}x 
    C(\sigma ^{j-2}) (x)   
     C(\int_{x}^{\infty} \sigma^{n-j+1} \sigma ' )( x) 
     dx  
   \ee
   \be
      + \frac{1}{2\pi^{2}} \int_{0}^{\infty}x 
    C(\sigma ^{n-j+1}) (x)   
     C(\int_{x}^{\infty} \sigma^{j-2} \sigma ' )( x) 
     dx  
   \label{eqn4.12}
  \ee
  or
   \be
      \frac{-1}{2 \pi^2} \frac{1}{n-j+2} \int_{0}^{\infty}x 
    C(\sigma ^{j-2}) (x)   
     C( \sigma^{n-j+2}( x) 
    \ dx  
   \label{eqn4.13}
  \ee
   \be
      +\frac{-1}{2 \pi^2} \frac{1}{j-1} \int_{0}^{\infty}x 
    C(\sigma ^{j-1}) (x)   
     C( \sigma^{n-j+1}( x) 
    \ dx  .
   \label{eqn4.14}
  \ee

  We can do the $j=2, j=n$ cases separately just as easily (the above formulas 
  are not even all required in that case) and putting the two cases together 
  and reindexing when necessary we arrive at the conclusion of the 
  theorem.
  
  Our final step is to extend this to functions other than powers. The 
  standard uniformity arguments used in the Wiener-Hopf theory apply here 
  if we can show that 
  $$||\tr f(B_{\alpha}(\sigma))-\tr B_{\alpha}(f(\sigma))||_{1}= O(1)$$ 
  uniformly for $\sigma$ replaced by $1-\lambda +\lambda \sigma$ and 
  $\lambda$ in some complex neighborhood of $[0,1].$ The details of this
  are found in \cite{W2}. The norm above is the trace norm. Given 
  sufficient analyticity conditions on $f$, it is only necessary to prove  
  $||B_{\alpha}(\sigma_{1}) B_{\alpha}(\sigma_{2}) -B_{\alpha}(\sigma_{1} 
  \sigma_{2})||_{1}=O(1) $ where the 
   $O(1)$ here depends on properites of $\sigma_{i}.$
  A trace norm of a product can always be estimated by the product of  
  two Hilbert-Schmidt norms and in this case we need to estimate the 
  Hilbert Schmidt norm of the operator with kernel
  $$X_{(1,\infty)}(z)\int_{0}^{\infty}\sigma_{i}(t/\alpha) \sqrt{xz}t 
  J_{\nu}(xt) J_{\nu}(tz) dt.$$ Using integration by parts, and  
  integration formulas for Bessel functions this is easily estimated to 
  be bounded. For analogous details see \cite{W2}. 
  Thus for suitably defined $f$ we can extend our 
  previous theorem to the more general case. The $f$ of interest is 
  $\log (1+z).$ This will satisfy the necessary analyticity conditions if 
  we consider small enough $k$. The necessary conditions are collected 
  in the following:
  
  \begin{theorem}
  Suppose $f$ is a real-valued function with $[\nu]+2$ derivatives all 
  contained in 
  $L_{1}$. Then for sufficiently small $k$ (say $k< 
  ||\sigma||_{\infty}^{-1}$) 
  $$\check{\phi}(k) \sim \exp \left\{ \frac{\alpha}{\pi}\int_{0}^{\infty}
   ikf(x)dx -\frac{ik\nu}{2}f(0) -\frac{k^{2}}{2 \pi^{2}}\int_{0}^{\infty}
   xC(f)^{2}(x)dx \right\}.$$
   \end{theorem}
   Proof: The form of the answer follows from the computation of the mean 
   given earlier and from the fact that the constant term in the previous 
   theorem is exactly half of the answer in Szeg\"{o}'s Theorem. Thus the 
   above answer for the $\log$ function must be half as well.
   
   The author would like to thank both Craig Tracy and Harold Widom 
   for many useful and helpful conversations.

\end{document}